\documentclass[12pt,notitlepage,twoside,a4paper]{amsart}
\usepackage{amsmath}
\usepackage{amsfonts}
\usepackage{amssymb}
\usepackage{amscd}
\usepackage{amsthm}
\usepackage{amsbsy}
\usepackage[french, british]{babel}

\usepackage{amsmath,amssymb,enumerate}

\usepackage{epsfig,fancyhdr,color}

\usepackage{amssymb}
\usepackage{amsmath,amsthm}
\usepackage{latexsym}
\usepackage{amscd}
\usepackage{psfrag}
\usepackage{graphicx}
\usepackage[latin1]{inputenc}
\usepackage[all]{xy}
\usepackage[mathcal]{eucal}

\definecolor{NoteColor}{rgb}{1,0,0}

\renewcommand{\textsc}{\textcolor{red}}

\newtheorem*{theorem 1}{\rm\bf Proposition 1}
\newtheorem*{theorem 2}{\rm\bf Proposition 2}

\theoremstyle{definition}

\theoremstyle{remark}

\def\interieur#1{\mathord{\mathop{\kern 0pt #1}\limits^\circ}}

\begin{document}
\title[Mathematics and  map drawing  in the eighteenth century]{Mathematics and  map drawing in the eighteenth century}

\author{Athanase Papadopoulos}
\address{Athanase Papadopoulos,  Universit{\'e} de Strasbourg and CNRS,
7 rue Ren\'e Descartes,
 67084 Strasbourg Cedex, France}
\email{papadop@math.unistra.fr}

\date{\today}
\maketitle

 
 \begin{abstract}

 We consider the mathematical theory of geographical maps, with an emphasis on the eighteenth century works of Euler, Lagrange and Delisle.   This period is characterized by the frequent use of maps that are no more obtained by the stereographic projection or its variations, but by much more general maps from the sphere to the plane. More especially, the characteristics of the desired geographical maps were formulated in terms of  an appropriate choice of the images of the parallels and meridians, and the mathematical properties required by the map concern the distortion of the maps restricted to these lines.
 The paper also contains some notes on the general use of mathematical methods in cartography in Greek Antiquity and in the modern period, and on the mutual influence of the two fields, mathematics and geography. The final version of this paper will appear in Ganita Bharati (Indian mathematics).

    \end{abstract}

\smaller{

\larger}

 \tableofcontents

   \section{Introduction}

\begin{figure}
\centering
 \includegraphics[width=10cm]{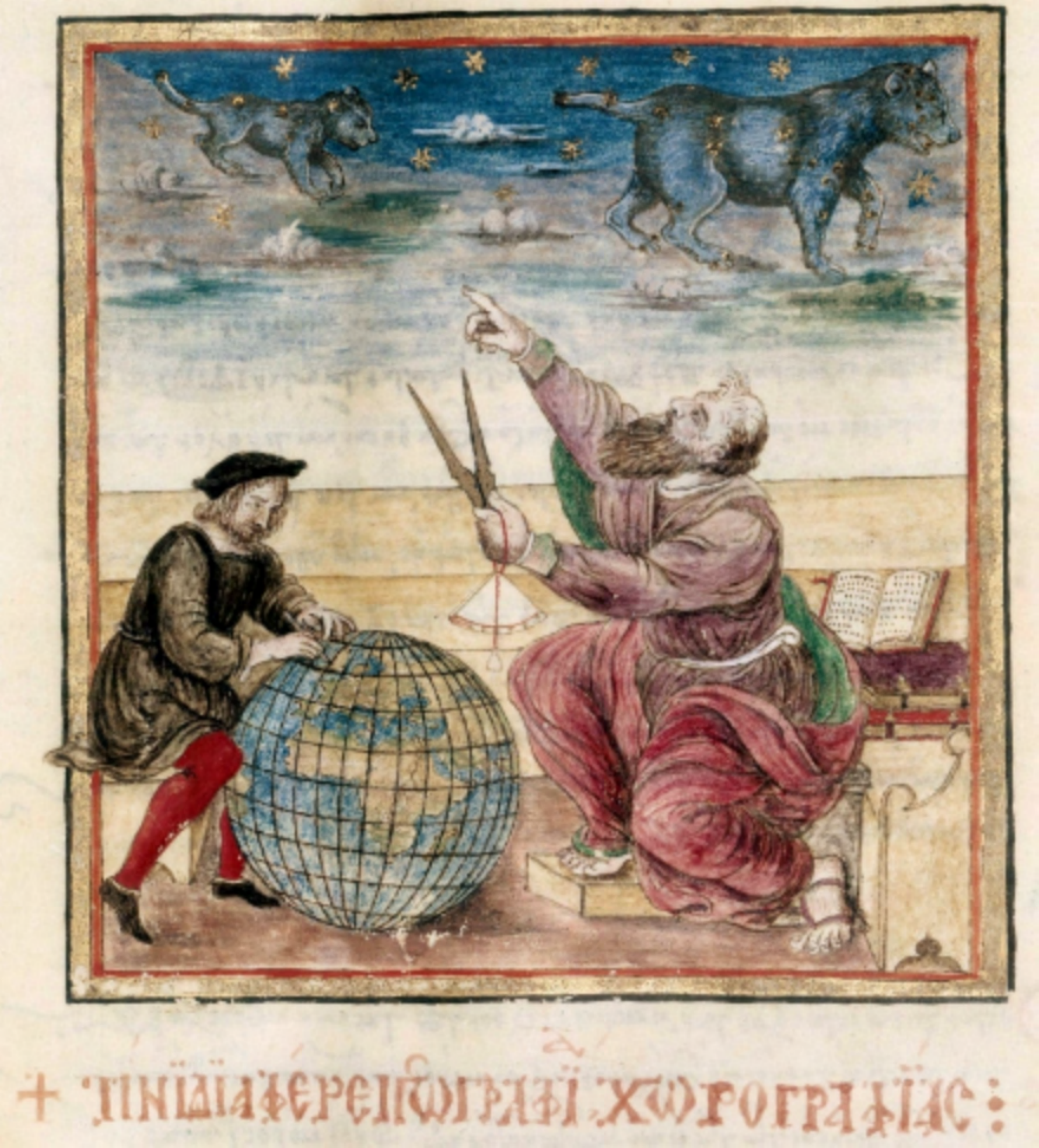}    \caption{\small The first page of a copy of a manuscript of Ptolemy's \emph{Cosmographia}.  Constantinople, end of XVth century. Biblioth\`eque Nationale de France. (Department of manuscripts, Codex Greek  No. 1401 fol. 2.).}   \label{fig:BN2} 
\end{figure}

  Cartography, or the art of drawing geographical maps, as it reached us, is rooted in the works of mathematicians of Greek Antiquity.
  It has close connections with geometry and astronomy.  
The Greek geographers of Antiquity that we shall mention considered the Earth as a spherical body, and they were able to carry out computations of large distances on its surface using methods of spherical geometry, in particular, of formulae making relations between the three sides of spherical triangles\footnote{We recall that a spherical triangle, unlike a Euclidean one, is completely determined by its three angles.  Thus, there are formulae giving the length of a side in terms of the angles. This fact was known to the Greek geometers. An important treatise on spherical triangles is Menelaus' \emph{Spherics} (1st-2nd c. BCE), see \cite{RR2}.} and between  angles and sides of such triangles, that is, spherical trigonometry (even though the expression was not used yet).  Spherical geometry was also used in astronomy, since the heavenly sphere was, like the Earth, assimilated to a sphere.\footnote{The radius of that sphere was irrelevant, and it could also be considered as infinite, since the distance between two stars was taken to be the angle they make from the observer's viewpoint.}   Furthermore, geographers from Antiquity realized that it was possible to measure large distances on the surface of the Earth by methods involving the study of the positions of stars and of distances between them. Thus, they developed methods for measuring   distances between two distant locations on the surface of the Earth using their knowledge in astronomy.  The illustration in Figure \ref{fig:BN2} is extracted from a 15th century byzantine manuscript of a copy of Ptolemy's  influential work, the \emph{Geography}. It represents the author, with the index of his right hand pointing to the Small Bear constellation, while measuring, with his left hand, angles between the stars, using a quadrant, a compass and a plumb-line, and dictating the corresponding distances on the surface of the Earth to another man, dressed as a medieval geographer, whose task is to mark the latitudes on a World map. On the visible face of this map are represented the three known continents, Europe, Asia and Africa. The two systems of lines on the sphere represented on this World map are the parallels and the meridians. They will be an important element of discussion in the present paper. The two grids that consist of a selection of (images of) meridians and longitudes that are represented on a geographical map are usually called \emph{graticules}. 

In this paper, instead of graticules we shall consider the complete systems of lines, and we shall take the liberty of using modern mathematical terminology, calling the parallels and the meridians \emph{foliations} of the sphere (with two singular points each, at the North and South poles). The images of these foliations by a geographical map are  foliations of the Euclidean domain which is the image of this map. 
 The leaves of the first foliation are the \emph{parallels} of the surface of the globe representing the  Earth, that is, they are small circles obtained as the intersection of the sphere with planes parallel to the one containing the equator. The latter is a great circle (the intersection of the sphere with a plane passing through its center) which is perpendicular to the rotation axis of the sphere; it constitutes the separation line between the so-called Northern and Southern hemispheres. The parallels are the equidistant lines to the equator.
The leaves of the second foliation are the \emph{meridians},  viz. the half-great circles perpendicular to the equator (or, equivalently, the half-great circles connecting the North and South poles). The foliations by meridians and parallels are perpendicular to each other. The images of these foliations by the geographical map may or may not be perpendicular, depending on the map. 

\begin{figure}[htbp]
\centering
\includegraphics[width=7cm]{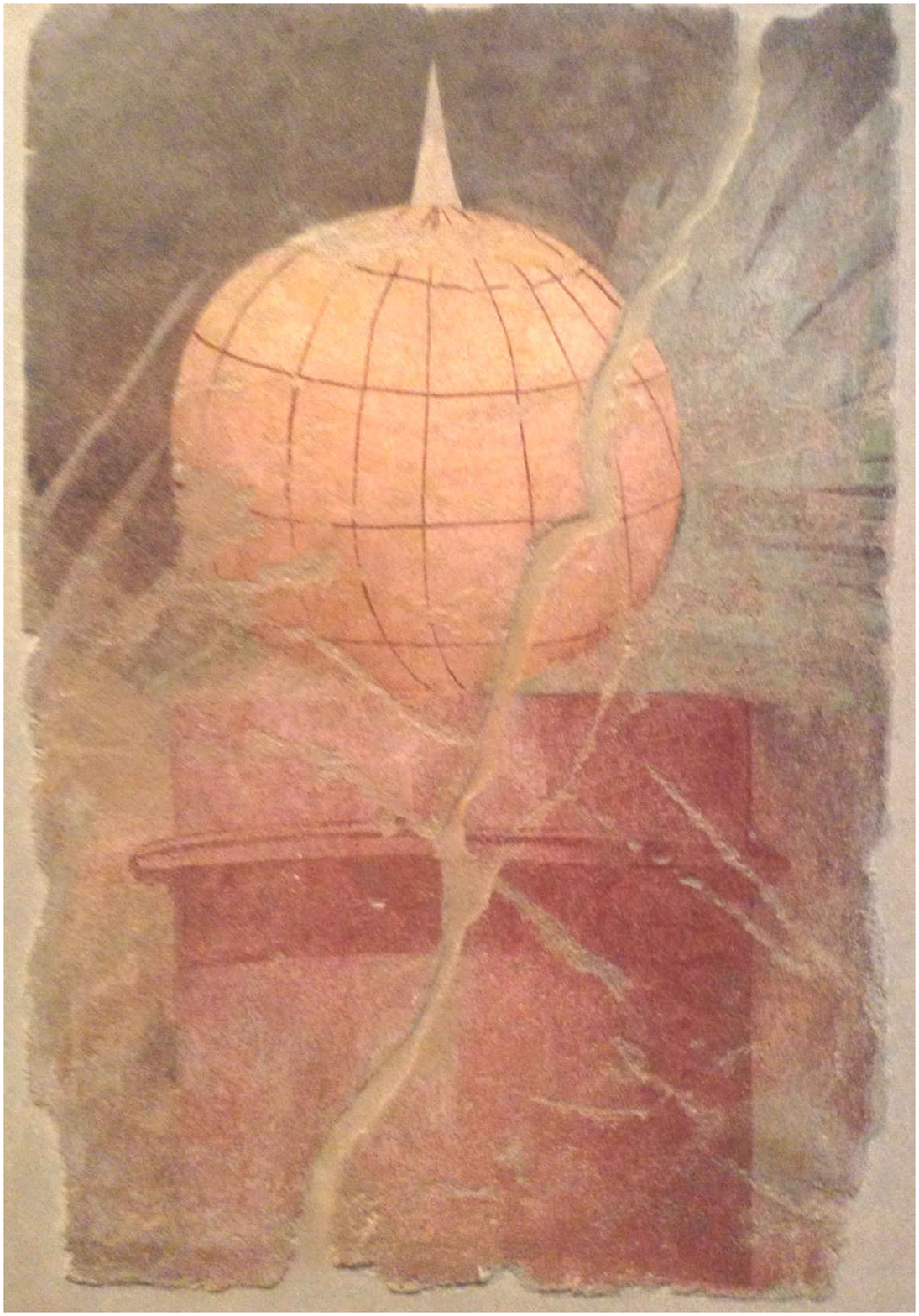}
\caption{\small Celestial globe with the two perpendicular foliations by parallels and meridians. Wall painting fragment from the peristyle of the Villa of P. Fannius Synistor at Boscoreale (Province of Naples), 
ca. 50-40 BCE. The fragment is exhibited at the Metropolitan Museum of Art, New York,  Department of Greek and Roman Art. (Photo A. Papadopoulos.)} \label{Met}
\end{figure}

Using these two systems of curves, the position of an arbitrary point on the surface of the Earth is usually determined by two coordinates: (1) the \emph{latitude}, which is the distance from this point to the equator, or, equivalently, the distance to the equator of the parallel on which the point lies; (2) the \emph{longitude}, which is the angle made by the meridian passing through the given point and a fixed reference meridian. Here, the angle made by two meridians is the dihedral angle made by the two planes containing them. Longitudes are taken between 0 and $\pi$, with the following convention: the longitude is called Eastern or Western according to whether the given point is at the East or West of the chosen reference meridian. 
 Ptolemy used extensively these coordinates in his \emph{Geography}. He took the reference meridian to be the one passing through his city, Alexandria, which in those times was a world major center of knowledge. Figure  \ref{Met} is a representation of a celestial globe dating from the 1st century BCE which is reduced to its most elementary elements: the foliations by parallels and meridians. The same picture could have been used for the representation of the Earth.

Let us stress on the fact that whereas meridians are geodesics on the sphere, parallels are not. The degrees of latitude are proportional to distances on the sphere, by a factor which does not depend on the chosen meridian, whereas the degrees of longitude are not proportional to distances.

    The geometers from ancient Greece knew that it is not possible to have a faithful representation of (part of) a sphere on a Euclidean plane. This follows immediately from properties of  spherical triangles, e.g. the fact that angle sum is not constant, or from other the properties that make these triangles very different from the Euclidean ones. One of the themes that I will stress in this paper is that drawing a geographical map while taking into account this impossibility 
was reduced to the question of drawing in an appropriate way the images of the foliations by parallels and by meridians.
This simple remark, which in a sense is an obvious one, turned out to be of paramount importance. It was formulated explicitly in the 18th century,  and it was at the center of the researches of Lambert, Euler, Lagrange and others.   We shall discuss this in the work of Lagrange, in \S \ref{s:Lagrange} below.
           
   The rest of this paper is divided into 6 sections, the main ones being Sections  \ref{s:Delisle} to \ref{s:Lagrange}. They concern the 18th century activity on mathematical map drawing. We shall review there the works of Delisle, Euler and Lagrange and a few others. Before that, in Section \ref{s:antiquity}, we shall make an overview of the work done before this period, mainly by Ptolemy, and we shall use ideas from there in our account of 18th-century cartography. The last section of this paper, Section \ref{s:Conclusion}, contains some perspectives and concluding remarks.
   
   While surveying the works of geographers and map drawers that are well-known mathematicians (Ptolemy, Euler, Lambert Lagrange, etc.), we  included a minimal amount of biographical elements concerning them. But we thought it was useful to include some more information on the life of Joseph-Nicolas Delisle, a famous astronomer and geographer at the Saint Petersburg Academy of Sciences  with whom Euler collaborated on cartography, assuming that his name is poorly known to mathematicians.
   
   An extended analysis of mathematical geography of the eighteenth century is made in the book \cite{CPG}.

In this paper, the translations from the French are mine.

     \section{Cartography in Antiquity}\label{s:antiquity}
  
  Cartography is a field of practical importance in which advanced mathematics was needed, and I will start by a few words of motivation on this subject.
  
   In Ancient times, the drawing of geographical maps, besides its theoretical importance for the advancement of knowledge, was motivated by practical questions. Maps were useful to travelers, explorers, traders, adventurers and navigators. The making of geographical maps was sponsored by political rulers who needed to have precise representations of the provinces they controlled. Maps were also useful to land owners, for tax and inheritance purposes.    Geographical maps were also used by conquerors and generals, who needed to have an image of the lands they intended to overmaster,
   and by others who were simply curious to know the scope of the known inhabited and uninhabited world.  For the cultured man, skimming geographical maps would show that the community or country to which he belonged was not the only one on Earth. We may recall here a dialogue by Claudius Aelianus (ca. 175-ca. 235) in which this Roman author relates a conversation between Socrates and the Athenian statesman Alcibiades: The philosopher, seeing the excessive pride that his  interlocutor took from the extent of the land he controlled,  showed him a World map and asked him to find the place of Attica. When Alcibiades found the place, Socrates continued: ``Now show me the lands you own." Alcibiades responded that they were not marked on the map. Socrates then told him: ``How can you be so proud of properties which do not represent even a point on Earth!"
\cite[Book II, chap. 28]{Aelianus}
    
   Among the many Greek mathematicians who were thoroughly involved in geography, and in particular in map drawing, we mention  Autolycus of Pitane  (4th c. BCE), whose work, \emph{On the moving sphere}, on astronomical phenomena, is the oldest Greek mathematical treatise which completely survives, Eratosthenes of Cyrene (3rd c. BCE),  whose name is associated with a measure of the length of an Earth meridian, Hipparchus  of Nicaea (2nd c. BCE), who is the first known mathematician who discovered formulae of spherical trigonometry and who used them for establishing astronomical tables, Posidonius of Apamea (2nd c. BCE), who carried out a measure of distances from the Earth to the Moon and other  large distances, and the great Claudius Ptolemy (3rd c. CE), whose \emph{Geography} \cite{Ptolemy-geo}, which we  mentioned above, marks the climax of the Greek contribution to the named field. In this work, Ptolemy, besides compiling the main achievements of his predecessors and correcting  them, used all the mathematical, optical and astronomical tools that were at his disposal in order to give the most possible precise indications that were useful for drawing geographical maps of the known world. The prominent 18th--19th century French scientist Pierre-Simon de Laplace, in a passage of his major work \emph{Exposition of the system of the world},  in which he explains the importance of geography in the sciences, attributes to Hipparchus the method of computing the coordinates of locations on the Earth by their latitude and longitude using astronomical data, namely, the eclipses of the moon. In the same passage, Laplace recalls that Hipparchus' works did not reach us because they were destroyed with the fatal burning of the  library of Alexandria, and that we know these works only through those of Ptolemy \cite[vol. 2, p. 505]{Laplace}.

Like his predecessors,  Ptolemy used, for his computations of distances on the Earth, theoretical methods involving a comparison between arcs of circles on the Celestial sphere with corresponding
arcs of circles on the surface of the Earth.  From the philosophical point of view, this possibility of computing the earthly distances from their counterparts on the Celestial sphere, which combines astronomy and geography into one single field of investigation, was a ground for Ptolemy's (and before him, of Hipparchus', Eratosthenes' and others') belief that the whole universe is one single body.

By the end of Book II of his major book on Astronomy, the \emph{Almagest}, after establishing his well-known tables of the coordinates of the celestial remarkable places, Ptolemy writes: ``This table of angles should be followed by the locations of the most famous cities, according to their longitudes and latitudes, computed according to the celestial phenomena that are observed in each of these cities. But we shall treat distinctly this interesting subject which belongs to geography [\ldots]  \cite[t. 1 p. 148]{Almagest-Halma}." Thus, Ptolemy kept his word, since he included his astronomical observations as a major tool in the computation the geographical coordinates that he published in his \emph{Geography}. 

    Ptolemy had at his disposal the formulae for spherical triangles that were discovered a few decades earlier by Menelaus of Alexandria. He was also aware of the mathematical properties of two types of projections: the one onto a plane passing through the center of the sphere and the one onto a plane tangent to the sphere, both centered at an intersection point between the sphere and a diameter perpendicular to that plane. 
    
  We already mentioned that in our discussion of geographical maps, an important factor will be the images of the parallels and the meridians. Let us note, as a preliminary observation, two special cases:
  Under a stereographic projection centered at the North  pole, the parallels are sent to concentric circles centered at the image of the South pole, while the meridians are sent to straight lines meeting at this image. Another projection,   called the \emph{gnomonic projection} and which was known since the times of Thales, has a similar property. This is a projection from the center of the sphere onto a plane tangent to the South pole. 
 Under this projection, parallels are sent to circles centered at the South pole and meridians are sent to straight lines  passing through the South pole.

     In the first book of the \emph{Geography}, Ptolemy discusses the mathematical methods that his predecessors used for drawing geographical maps and explains how these methods can be improved.  Half of his treatise (in number of pages) is dedicated to tables of coordinates of known places; he indicates there about 8000 names of locations, from Ireland in the West to China in the East, each  with its coordinates, expressed in degrees and minutes of latitudes and longitudes. This work also includes descriptions of several big rivers, with the location of their sources and a specification of their form (curvature, etc.). Regarding India, Ptolemy indicates 270 locations, giving also the duration of the longest day for 30 cities in this country. The last book of the \emph{Geography} contains instructions for map drawing  (10 maps for Europe, 4 for Africa and 12 for Asia).

   During the Middle Ages, Ptolemy's  \emph{Geography} was translated  into Arabic by the mathematician  al-Khw\=arizm\=\i , and later, into Latin, from the Arabic. One may also mention here that Ptolemy wrote another geographical work, known in Latin as the  \emph{Planisphaerum}, which reached us only through Arabic translations (the Latin version as well as all the subsequent versions were translated from the Arabic;  the Greek original is lost).

\begin{figure}
\centering
 \includegraphics[width=\linewidth]{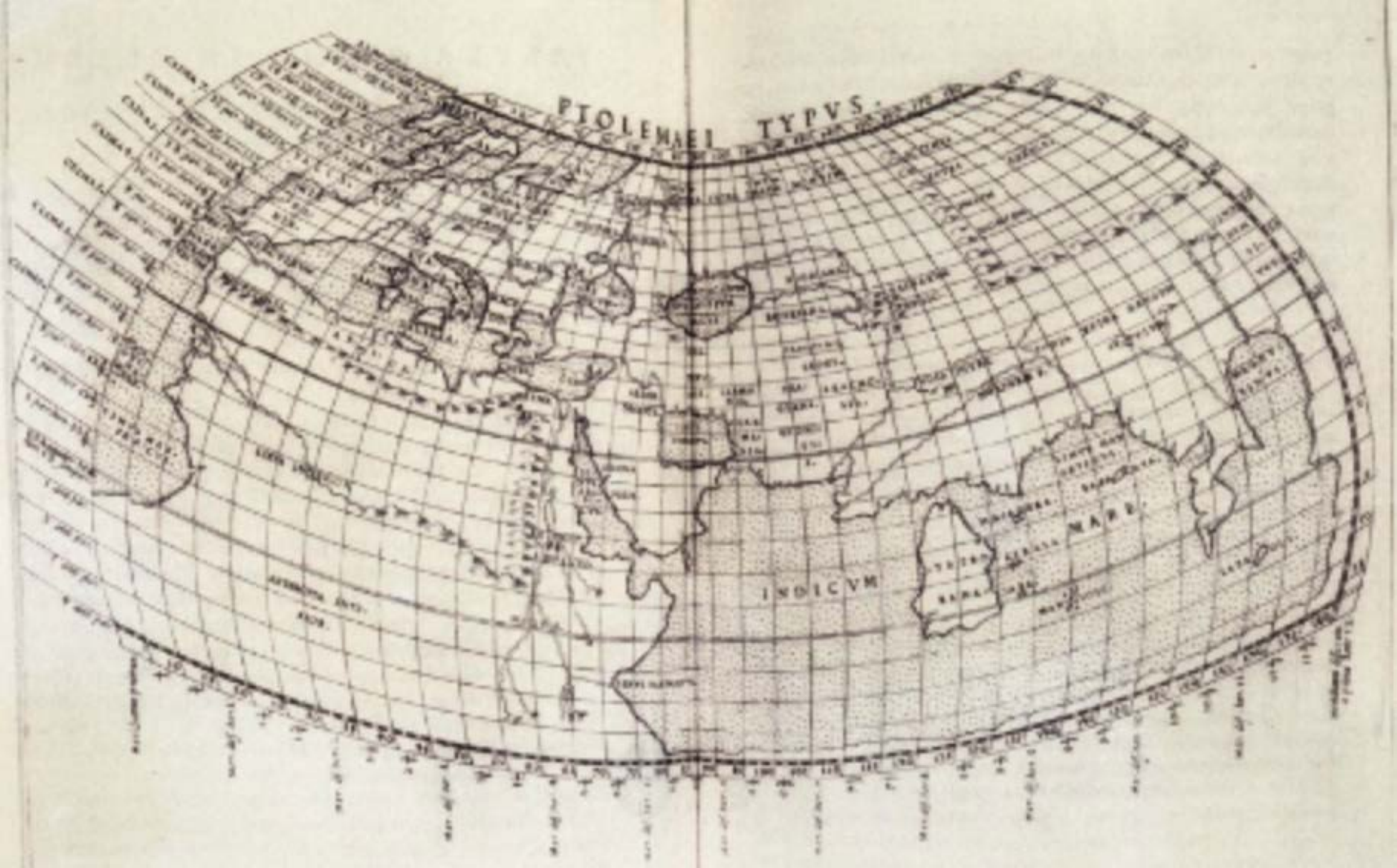}    \caption{\small Map of the Ecumene, or the ``Known World". From a Latin translation of Ptolemy's \emph{Geography}. Translator Girolamo Ruscello. Venice, around, 1564. Bodleian Libray, University of Oxford, Collection of Medieval and Renaissance Manuscript Illumination, Byw. H 5.9}   \label{fig:Ptolemy3} 
\end{figure}

Ptolemy's geographical maps did not survive but they were reconstructed from the coordinates he gave, first by Byzantine monks at the end the 13th century and later by geographers and artists, when his works reached the Latin world.  The first printed Greek version of the  \emph{Geography} was edited by Erasmus in 1553.
     Figure \ref{fig:Ptolemy3} represents a map of the Known World (``Ecumene"), from a 16th-century  edition of this work. One may notice that in this representation, the parallels and the meridians are perpendicular, and that the parallels are circles whereas the meridians are not.  We shall thoroughly comment on such properties in the rest of the present paper.

\begin{figure}
\centering
 \includegraphics[width=\linewidth]{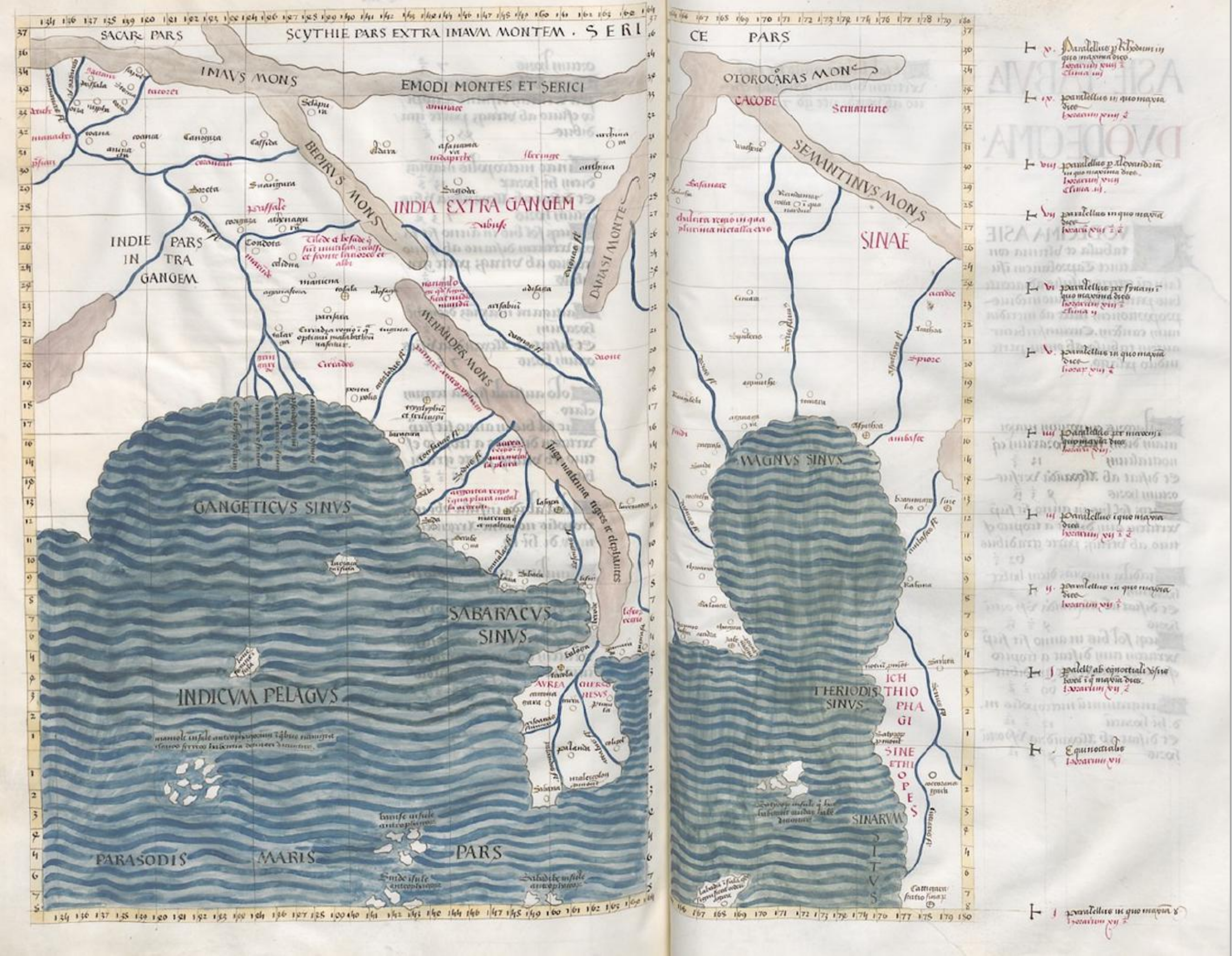}    \caption{\small Map of the Northern part of India, from a Latin edition of Ptolemy's \emph{Geography}. Translator Jacobus Angelus, around 1485.     Biblioth\`eque Nationale de France, Department of manuscripts. Latin No. 4804.}   \label{fig:Ptolemy-India}
\end{figure}

Figure \ref{fig:Ptolemy-India}, reproduced from a 15th century edition of the \emph{Geography}, represents a map of the Northern part of India. It carries the indication \emph{India extra Gangem} (India beyond the Ganges). The region depicted is bordered on the West by the Ganges Delta and on the East by the land of \emph{Sinae} (China). In this map, the meridians and the parallels are straight lines. The map,  from a 15th century Latin manuscript of the \emph{Geography}, was drawn according to Ptolemy's indications. On the right hand side of the picture are written the names of the parallels and their latitudes.

It was well-known, since mathematicians understood the first properties of spherical triangles, that it is not possible to draw maps from (a region of) the sphere to the Euclidean plane that preserve distances up to scale.
In general, (images of) parallels and meridians are highlighted as significant networks of lines on the drawing, and it was natural that geographers seek for maps  in which distances are preserved along these networks, or at least along one of the them. In the next sections, we shall review the way in which the 18th-century mathematicians-geographers dealt with such problems.   

\begin{figure}[htbp]

\centering
\includegraphics[width=12cm]{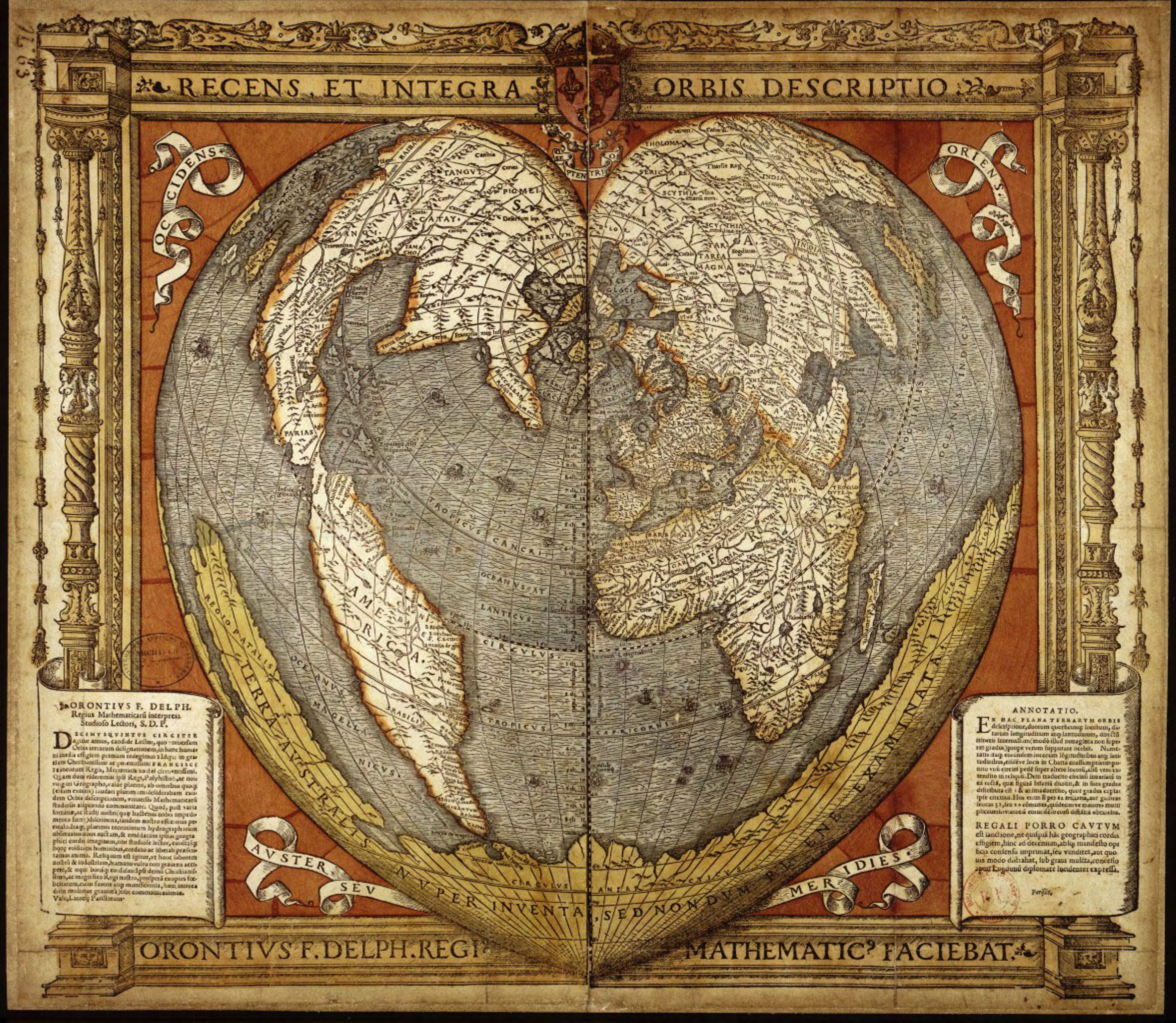}
 \caption{\small A world map, drawn between 1534 and 1536 by Oronce Fine, Biblioth\`eque Nationale de France.} \label{Coeur}
\end{figure}

 To end this section, let us look at another map on which the foliations of parallels and meridians are represented in a peculiar way. This is the 16th-century heart-shaped map represented in   
  Figure \ref{Coeur}, carrying the name \emph{Recens et integra orbis descriptio} (Recent and complete description of the world), drawn by the French mathematician and astronomer Oronce Fine  (1494--1555). The latter  was the first to hold the chair of mathematics at the Coll\`ege Royal (which became later the Coll\`ege de France). In this map, the parallels constitute a foliation whose lines are close to circles near the North Pole and to straight lines near the South Pole. The meridians are everywhere almost perpendicular to the images of the parallels.   Distances along the parallels and along the central meridian are preserved up to scaling. Such a projection is also known under the name \emph{Stabius--Werner projection}, in honor of two cartographers, Johannes Stabius (1460--1522), who was the first to highlight it, and Johannes Werner (1466--1528), who wrote a treatise on it.   Besides the remarkable properties that this geographical map satisfies on the meridians and parallels, its shape  has a practical advantage, namely, it leaves in the dark the question of whether North America and Asia are connected or not, a question whose response was unknown at that time. The drawing of this map is also based on Ptolemy's tables,  and its conception uses advanced spherical geometry and geometrical constructions combining the stereographic and the gnomonic projections; see the discussion in \cite{Hapgood}.

We next pass to  eighteenth century geography.

   \section{Joseph-Nicolas Delisle} \label{s:Delisle}
  
   We shall review in some detail the work of Joseph-Nicolas Delisle, and we start with a few words on his life.
   
   Born in Paris in 1688, Delisle was admitted in 1714 at the Royal Academy of Sciences  as a training astronomer. Soon afterwards, he published several memoirs on astronomy and physics.
 In 1725, the tsar Peter the Great, who was aware of Delisle's talents, invited him officially to Saint Petersburg, proposing him the post of director of the astronomical department at the Academy of Sciences that he was planning to found in that city. Delisle accepted and, a few months later, moved to Saint Petersburg. Louis the 14$^{\mathrm{th}}$,  King of France, who died the same year, promised to Delisle that his situation at the Royal Academy, as well as a position he had at the Coll\`ege Royal, on a Chair  of Mathematics, will be kept for him until his return.

Delisle  founded an astronomical school in Saint Petersburg which became, a few decades later, one of the most renowned in the world. He also supervised the construction of an astronomical observatory of which he was appointed the director, on the Vasilievsky island of Saint Petersburg. The observatory is situated on the top three floors of the Kunstkamera. It became soon one of the most famous in Europe.   At the same time, Delisle was in charge of the department of geography at the Academy. The principal task at this department was to make precise measurements of the territories and to draw maps of the  huge Russian Empire.

In 1736--1737,  Delisle proposed, for the elaboration of a map of the Russian Empire, a general triangulation of these lands. Triangulations, for measuring distances on the Celestial sphere, were already used in Antiquity by Euclid, Heron and others, see  e.g. \cite[p. 845, 893.]{Neu-H}. General maps of France, the Netherlands and other European countries had already been drawn using triangulations. This  was the first time that Russia was submitted to precise measurements using this method. 

For that purpose, Delisle organized several expeditions to the various parts of the Empire, from Central Asia to the remote frozen lands. Distances between principal locations needed to be measured based on astronomical observations. 
 In his computations, Delisle took into account the fact that the  shape of the Earth is  closer to a spheroid (a surface obtained by rotating an ellipse around an axis) than to a sphere. As a matter of fact, he considered that the Earth is flattened at the Poles, following a hypothesis emitted by Newton, which was in opposition to the theory emitted by the Parisian geographers and astronomers under the leadership of the Cassinis who thought that, on the contrary, the Earth is flattened at the Equator. 
 
 In 1737, Delisle read a paper on the triangulation of the Russian Empire to the Saint Petersburg Academy of Sciences. The paper was translated into English and published in the  Philosophical Transactions of the Royal Society of London \cite{Delisle-Proposal}. A short summary of this paper, highlighting a few passages, may be useful to our subject:

In the introduction, Delisle expands on the necessity and exigencies of geography, the drawing  of charts and the measure of distances on the surface of the Earth, and he recalls how this was dealt with at various stages of history, and how the need for more and more precise measurements arose. He emphasizes the ``immense labours of modern mathematicians" in this field, mentioning the works of  the mathematicians Fernel and Snellius, of the astronomer Riccioli and others, who computed lengths of meridians based on the laws of geometry and using astronomical observations. He stresses on the importance of conducting such measurements at several locations  and comparing them. He recalls that there were contradictory conjectures on the exact shape of the Earth: some scientists considered it as flattened at the poles, and others at the equator. He gives some details on the units of measurements used in Russia, comparing them with the French and the English ones, explaining that what remains to do towards the perfection of geography in Russia is to employ these units in actual measurements, ``and to construct the charts by the most exact methods of geometry, taking care to set them down right, as to their true bearings, and to regulate them by the most exact astronomical observations of longitude and latitude that can possibly be made." He then declares that it will not be possible to reach this goal unless an equal and even greater work than that which has been undertaken in France and elsewhere is undertaken, towards the measurement of the Earth. He explains that this will be a long and painstaking process. Setting in a proportionate way the degrees on the meridians and on the longitudes should take into account the fact that the Earth is not perfectly spherical. He writes: ``In all this there might possibly arise an infinite variety, according to the figure the Earth might have; and as it is not yet decided what is the Earth's true figure, and that there is no better method of ascertaining it than by observations made in so great an extent as that of Russia." He then refers to the various experiments and the work done by Newton who published them in his 1686 edition of the \emph{Principia}, showing that the Earth is flattened at the poles. He recalls that  Huygens supported Newton but that Picard, Snellius, Eisenschmid, the Cassinis (Giovanni-Domenico and his son Jacques) and others refuted Newton's claim and carried out measurements that were supposed to prove that, on the contrary, the Earth is longer at the poles. He also mentions the mathematician and astronomer Dortous de Mairan who, in 1720, tried to reconcile the two points of view. 
 
 Delisle then recalls that back in 1720, while he was still in France, he imagined a new method of tackling this question, by observing the degrees of the parallel compared to those of the meridian.
He explains his method of measurement, and he points out the errors of the French, due to the fact that the observations they made were not sufficient. This method consisted in forming triangles based along the parallel of Paris and observing at the two ends the difference of the meridians. He declares that at that time he was not able to complete his design ``for want of alliance and for reasons which [he] shall pass in silence." 

Delisle then explains the method of constructing maps of a large country using triangles, in different directions,  ``linked together by means of objects  seen successively one from another", and he  describes the measurements that Jacques Cassini undertook in 1733, which led to contradictory results, because they were based on some error in the initial information on longitudes, due to some mistakes in the observation of Jupiter's false satellite. He recalls that it became clear that to draw correct maps of France, it was necessary to make observations in regions far from France, and that this gave rise to the French expeditions to Peru, Lapland and other places, involving astronomers and mathematicians. At the moment where Delisle was writing his report, he was still not aware of the exact results of these expeditions.

Delisle then declares that in view of all this, he thought that it was necessary to undertake a work of the same nature in Russia, 
with an advantage that the other nations do not have, namely, the vast extent of the Russian Empire in every direction.  He then expands on where the vertices of the triangles should be taken: generally, in open places where the ice is even or where there are not trees, so that light signals that are used for indicating the locations can be perceived from far away. He concludes his memoir by describing several  astronomical instruments that are needed for this task, taking into account whether the sides of the triangles used in the measurements terminate at places that are high or not with respect to the level of the sea.

%
%

Delisle belonged to both departments at the Academy of Sciences, astronomy and geography.  These departments  were part of  the ``class" of mathematical Sciences, and Delisle had close relations with mathematicians.
 During his first years in Saint Petersburg,  he worked with Jakob Hermann, who was, like him, a founding member of the Academy. Hermann was a mathematician who, like Euler, came from Basel;  incidentally, he was a relative of the Euler's mother. In 1731, Hermann, who was 29 years older than Euler, returned to Basel, and  Delisle started to search for another mathematician who would assist him in his work, and the natural choice was Euler.
 
%
 
When Euler and Delisle started their collaboration, in 1735, Euler was already involved in astronomy. In 1732, he had presented to the Academy a memoir titled \emph{Solutio problematis astronomici ex datis tribus stellae fixae altitudinibus et temporum differentiis invenire elevationem poli et declinationem stellae} 
(Solution to problems of astronomy: given the altitudes and time differences for three fixed stars, to find the elevation of the pole and the declination of the stars)  \cite{E14}.\footnote{The memoir was published in 1735. Let us note that, in general, there was a (sometimes very long) lapse of time of several years between the time when Euler presented his memoirs and the time where they appeared in print. This was essentially due to the large backlog that the journal of the Academy accumulated, due to the amount of writings they received from Euler, who was extremely prolific.} In this memoir, Euler used spherical trigonometry in order to solve the problem announced in the title; in fact, the article starts with the spherical cosine formula in a triangle $ABC$:
    \[\cos A= \frac{\cos BC-\cos AB\cos AC}{\sin AB\sin AC}.\]
    
      In 1735, Euler presented another memoir on astronomy,  titled \emph{De motu planetarum et orbitarum determinatione}
(On the motion of planets and orbits) \cite{E37}.  His collaboration with Delisle boosted his interest in the field, and it was under the latter's influence that he wrote his first book on astronomy,  \emph{Theoria motuum planetarum et cometarum}
(Theory of the motions of planets and comets) \cite{E66}
 and his first important memoir on the motion of the moon, \emph{Theoria motus lunae exhibens omnes eius inaequalitates}
(Theory of the motion of the moon which exhibits all its irregularities)
 \cite{E187}.  Euler worked on topics related to the subjects of these two memoirs until the end of his life. 
  We shall talk more about Euler's work on geography in \S \ref{s:Euler} below.

The year 1740 saw the beginning of a conflict between 
Delisle and the Academy's administration. One reason was related to the so-called \emph{Atlas Russicus} (the ``Russian Atlas"), a project which was initiated by Peter the Great, and of which Delisle was in charge. Delisle kept postponing the publication of this atlas, and in 1740 this charge was taken away from him and given to Euler, who became in effect the head of the department of geography. To situate the problem in its context, one should remember the international competition that was taking place in this period, for drawing and publishing maps of Russia and the bordering countries (in particular China and Japan). The monarchs, who were directly supervising the Academy, were eager to see the atlas published.

Euler spent one year working on the atlas, in collaboration with the German mathematician, geographer and astronomer Gottfried Heinsius who was settled in Saint Petersburg, and without Delisle, until he moved to Berlin, in June 1741. During his long stay in Berlin and until Delisle returned to France, Euler kept informing the latter of his astronomical discoveries, and in particular his researches on the moon. The correspondence between Euler and Delisle during Euler's stay in Berlin shows that each of the two men had a high respect for each other's work. Euler, from Berlin, continued to follow the evolution of the atlas, whose direction was given to Heinsius. The Atlas was eventually published in 1745. It consists of 20 maps.  
    
 Later, the situation between Euler and Delisle  deteriorated. In a letter to Johann Kaspar Wettstein,\footnote{Johann Kaspar Wettstein (1695-1759) was, like Euler, from Basel, and the two men were friends when they were young.  Later, Wettstein moved to England where he became chaplain to the royal family. Euler and Wettstein kept close relations.} dated 5 June 1753 \cite[p. 443]{Opera-4-A-7}, Euler relates the criticism directed by Delisle towards the project the Russian Atlas that he was leading, writing that the latter was treating the project with despise. In the same letter, Euler tells Wettstein that  the Russians were unhappy with Delisle's memoir on the far-Eastern regions of Russia (``beyond Kamchatka"), especially for what concerns lies and voluntary errors contained in the latter's memoir on the Kamchatka expedition (see the notes in \cite[p. 445 and 457]{Opera-4-A-7}).

A short biography of Joseph-Nicolas Delisle is recorded in the memorial article \cite{Fouchy} by Jean-Paul Grandjean de Fouchy, astronomer and  perpetual secretary of the Paris Academy of Sciences, and as such, in charge of writing the necrological articles about the members of this Academy. One may find other information on Delisle  in his correspondence with Euler and with various scientists.

\section{Euler's writings on Delisle's method}
\begin{figure}
\centering
 \includegraphics[width=\linewidth]{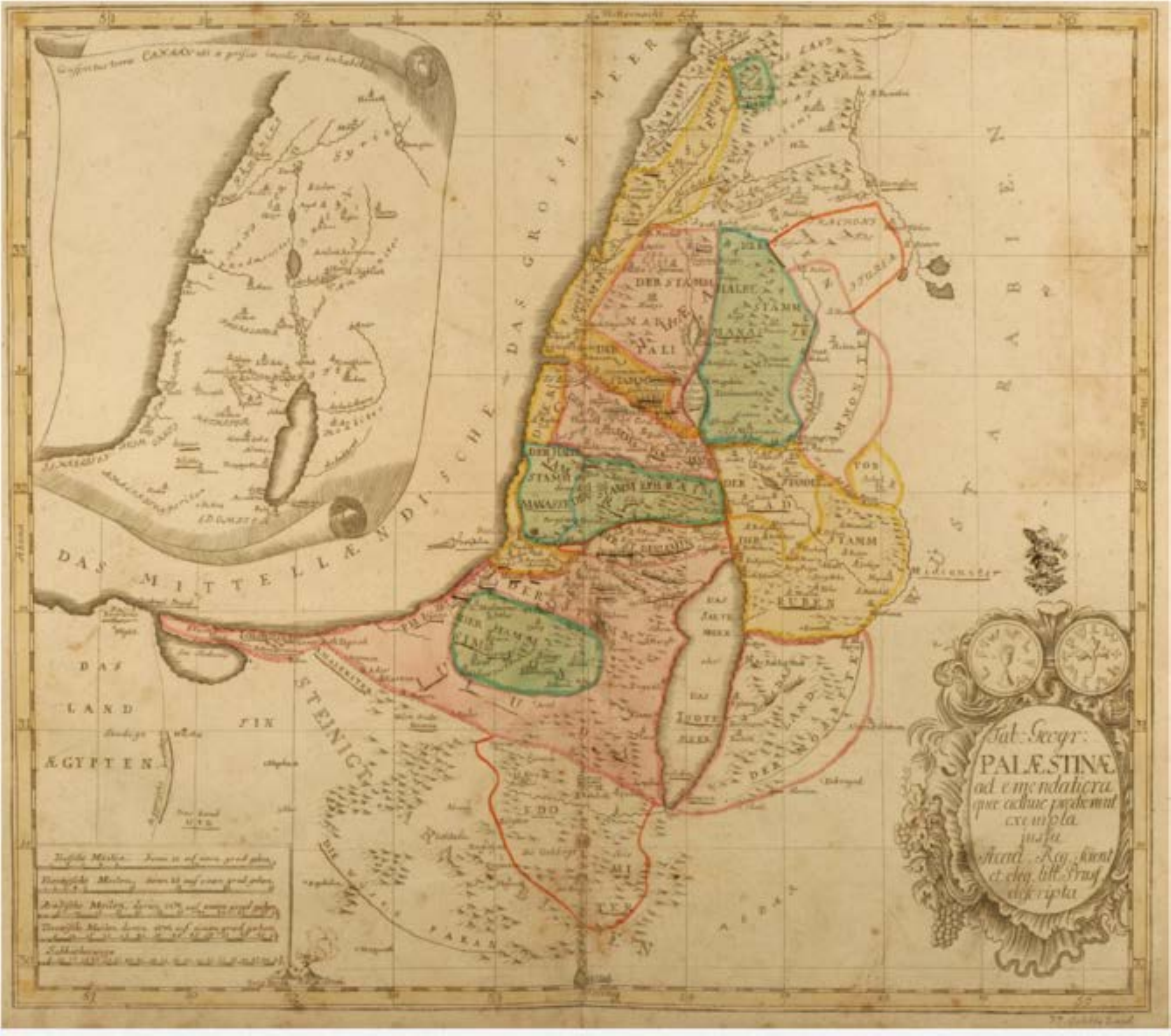}    \caption{\small Map of Palestine, from Euler's \emph{Atlas Geographicus} (Berlin, 1753)}   \label{fig:Euler-Palestine}
\end{figure}

\begin{figure}
\centering
 \includegraphics[width=\linewidth]{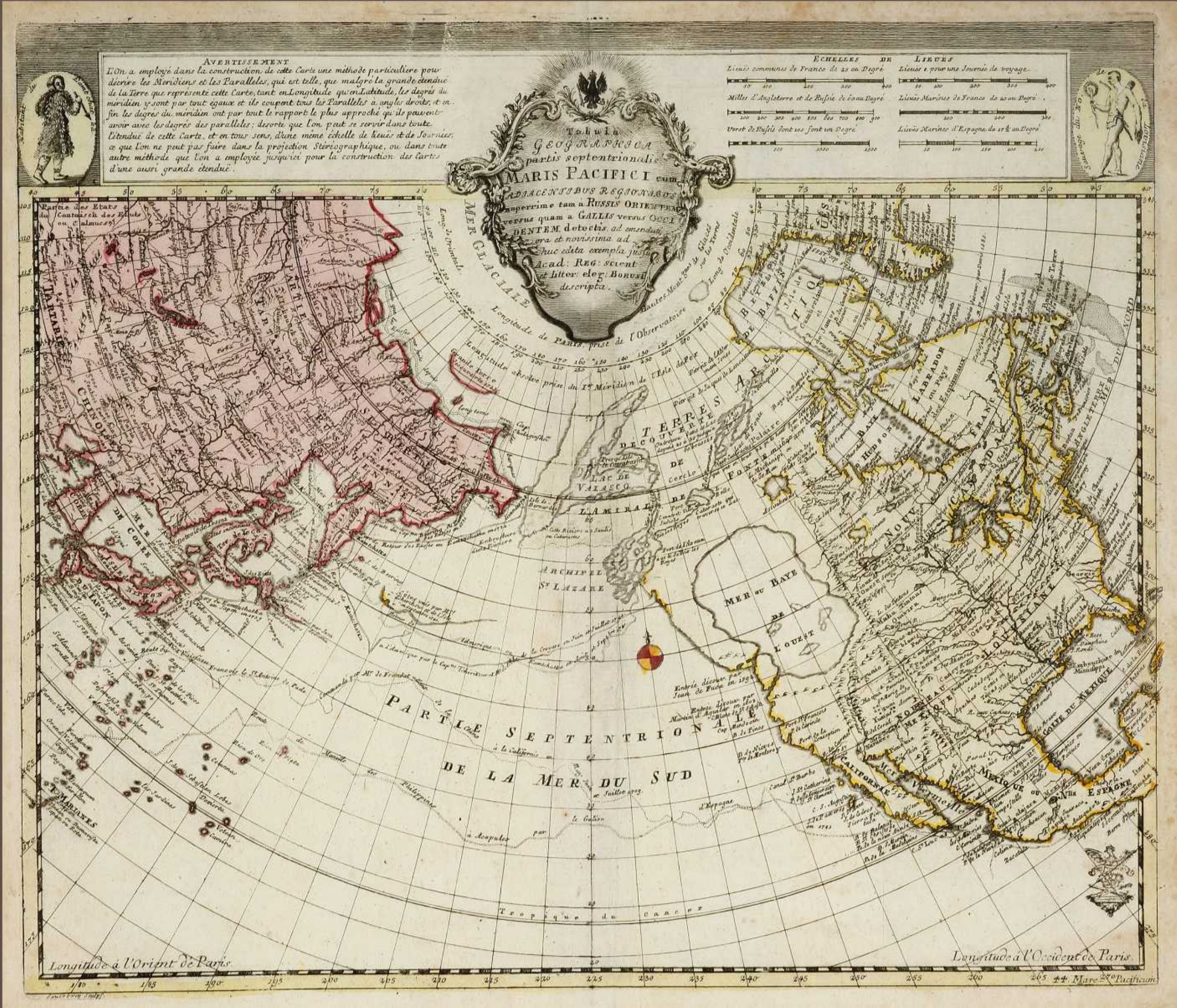}    \caption{\small Map of the Northern Pacific, with the Eastern part of Asia and the Northern part of America, from Euler's \emph{Atlas Geographicus} (Berlin, 1753)}   \label{fig:Euler-Map41} 
\end{figure}

 In this section, we shall review two texts by Euler on Delisle's method of map drawing. The first text is extracted from the  \emph{Atlas Geographicus omnes orbis terrarum regiones in XLIV tabulis exhibens} (Geographical atlas representing in  44 maps all the regions of the Earth), published by the \emph{Acad\'emie Royale des Sciences et Belles Lettres de Prusse}, in Berlin, in the year 1753  \cite{Atlas}. This atlas contains 44 maps, and it was edited under the direction of Euler who also wrote the preface. Several projections are used in drawing the maps of this atlas. It concerns  only large parts of the Earth globe: besides maps of the whole known world, it contains maps of regions of the size of a country. All the projections used in drawing these maps satisfy  the following two properties: the meridians are perpendicular, or almost perpendicular, to the parallels, and the degrees on longitudes are proportional to the degrees of latitudes. It is mentioned on the cover page that this atlas is principally meant for its usage in schools. We have reproduced, in Figures \ref{fig:Euler-Palestine} and
 \ref{fig:Euler-Map41},  two maps from this atlas, one of Palestine, in which the meridians and the parallels are straight lines, and one of the Northern part of the Pacific, together with the regions from  Asia (Eastern part of Russia) and America that surround it, in which the parallels are circles and the meridians are straight lines.

 We translate part of the preface of the atlas, which concerns the map of Figure  \ref{fig:Euler-Map41},  the last one in the series, which is drawn using Delisle's method. Euler, in this passage, comments on  this method. The preface is dated Berlin May 13, 1753. Euler writes:
\begin{quote}\small 
In the description of this map, we have kept the method which the famous  Mr. Delisle used, which seems to us the most appropriate for a good representation of these Northern regions. We have used the same method as in the general map of the Russian Empire published by the Saint Petersburg Academy, which at first glance may not be satisfying; we shall explain it here in a few words.

In this representation, all the meridians make straight lines and all their degrees are equal: two meridians that are distant apart by one degree converge in such a manner that under two elevations from the pole the degrees of longitude make with the degrees of latitude the same ratio as in reality. For this, the two elevations of the pole that we choose are those that are at the same distance from the extremities of the region that we wish to represent, than from its middle. In this manner, it happens that under these elevations the ratio between the degrees of longitude and latitude happens to be accurate, and that at the other locations it is almost correct. And thus, the whole representation makes the positions as much precise as possible. In the map of the Russian Empire, and in the present one which represents the regions situated between Asia and North America, we have chosen the elevations of 60 and 45 degrees from the Pole, which were the best suited for that purpose. In reality, all the meridians merge at a point, but this point is not the Pole; it is 7 degrees farther than the pole, and from this point as center all the parallel circles are described with equal distance among themselves. Now if in a similar frame we mark all the regions contained between the elevations of 45 to 68 degrees from the pole, their positions will differ by such a small amount from the true ones that we can hardly see the error. But if we wanted to place in this frame regions that are closer to the Pole or to the Equator, the error would be enormous. Hence, we see that we must completely exclude the representation of these regions, as Mr. Delisle did in the Map.

One should not regard as a shortage of this map the fact that the center in which all the meridians intersect is so far from the Pole, since according to this method it cannot be marked. Also, we should not be surprised if on this map the parallels, which form semi-circles, do not occupy 180$^{\mathrm{o}}$ in longitude, but much more, sometimes even up to 250$^{\mathrm{o}}$. From this we see that this map does not suffer from the fact that it has a too large amount of longitude defect. 
 \end{quote}

 Another commentary on Delisle's method, and in particular on the fact that the meridians do not meet at point which corresponds to the North or South pole, is given by Euler in a memoir we shall review now. 

In the year 1766, Euler, who was 69 years old and in his highest production period, presented to the Saint Petersburg Academy of Sciences three  memoirs on geography (the memoirs were published a couple of years later, \cite{Euler-pro-Desli-1777, Euler-rep-1777, Euler-pro-1777}). We review here one of them, which concerns Delisle's method of projection. The memoir is  titled
 \emph{De proiectione geographica Deslisliana in mappa generali imperii russici usitata} (On Delisle's geographical projection used in the general map of the Russian empire)  \cite{Euler-pro-Desli-1777}. We shall talk about the two others in the next section.

Euler starts his memoir by stating the inconveniences of a stereographic projection that was used by Hasius,\footnote{Johann Matthias Hasius (or Haas) (1684-1742) was a professor of mathematics in Wittenberg. He drew several maps and published a treatise on cartography, \emph{Sciagraphia integri tractatus de constructione mapparum omnis generis} (M\"uller, Leipzig, 1717). We mention his map  ``Imperii Russici et Tatariae universae tam majoris et asiaticae quam minoris et europaeae tabula" (Geographical map of the Russian Empire and of Tataria, both large and small, in Europe and Asia), published in Nuremberg in 1739, which is related to our topic here.} also known as the ``stereographic horizontal projection", or the ``Haas stereographic projection". This is a projection whose center is a point on the sphere, onto a plane which intersects the sphere along a small circle. The point diametrically opposite to   the center of the projection is called the ``location", and the projection plane intersects the sphere in a circle called the ``rational horizon"  of the location.  The meridian passing through the location is called the ``central meridian". In \S 2 and 3 of the memoir, Euler mentions two inconveniences of this map: the excessive unevenness of the scale on the image of the central meridian, and the fact that the images of the meridians are not evenly curved on the map. He then writes  (\S 3) that a certain property of a geographical map is desirable, namely, that one should be able to extract any part of it in such a way that this part is a reasonable map of the subregion represented, and can be used without change. He then says that this implies that the difference in curvature of the meridians should not be noticeable. 
The stereographic projection, in which all meridians are represented by straight lines that intersect at the pole, does not have the latter inconvenience, but it has the shortage that the scale is very uneven along these meridians.

In \S 5, Euler states four properties that are ideally required from a geographical map:
\begin{enumerate}

\item The images of the meridians are straight lines; 
 
\item the degrees of latitudes do not change along meridians;
 
\item the images of the parallels meet the images of the meridians at right angles; 

\item at each point of the map, the ratio of the degree on the parallel to the degree on the meridian is the same as on the sphere.
\end{enumerate}

He then declares that since this cannot be achieved simultaneously, one may request,  instead of the last condition,  that the deviation of the degree of latitude to the degree of longitude at each point from the true ratio be
 as small as possible (ideally, this error should be unnoticeable).

After this, Euler recalls that Delisle was in charge of constructing such maps, and that in doing this he started by requiring that the ratio between latitude and longitude be exact at two noteworthy parallels.
He writes that Delisle was of the opinion that the deviation will be small everywhere if these two parallels are equidistant from the middle parallel of the map and from its outermost edges. He then says that the question becomes that of choosing these two circles of parallels in such a way that the maximum deviation over the entire map is minimized.

Starting from \S 7,  Euler explains, based on mathematical considerations, how to construct a family of  straight lines representing meridians which are one degree distant from each other. The construction starts with the choice of a meridian passing through the Russian Empire, considered  as the principal meridian. He then gives a method for drawing the other meridians, which become straight lines that all intersect at a point which  is not the North Pole of the Earth. This implies that near the North Pole, the picture is wrong, but Euler notes that for the drawing of maps of the Russian Empire, there is no need to represent regions beyond 70$^{\mathrm{o}}$ of latitude. He computes the error in scale at this degree of latitude and he shows that it can be neglected.

Once the intersection point of the meridians is found, the parallels are drawn as circles centered at this point. The construction allows the construction of circles such that for the  region comprised between the two latitudes we started with, the ratio of the degree of latitude to  the degree of longitude is faithfully represented. 

In \S 10--16, Euler gives the mathematical details showing how far this representation differs from reality at the extreme points he started with, and in \S 17--23, he makes the actual computations in the  special case where the map is that of the Russian Empire.

In \S 24 and 25, he studies images of great circles on the sphere by the geographical map. He shows that these images are not noticeably different from Euclidean circular arcs and he computes the radius of such an arc, which he finds very large. Since a Euclidean circle of large radius is approximately a Euclidean line, he concludes that the shortest lines on the map do not differ sensibly from straight lines.

    \section{Euler as a geographer} \label{s:Euler}
    
Among Euler's early works on geography, we mention his memoir \emph{Methodus viri celeberrimi Leonhardi Euleri determinandi gradus meridiani pariter ac paralleli telluris, secundum mensuram a celeb. de Maupertuis cum sociis institutam}
(Method of the celebrated Leonhard Euler for determining a degree of the meridian, as well as of a parallel of the Earth, based on the measurement undertaken by the celebrated de  Maupertuis and his colleagues) \cite{E132},  presented to the Academy of Sciences of Saint-Petersburg in 1741 and published in 1750. 

From a more theoretical point of view, Euler published in 1777 three memoirs on mappings from the sphere onto the Euclidean plane, motivated by the problems of cartography:
\begin{enumerate}
\item  \emph{De proiectione geographica Deslisliana in mappa generali imperii russici usitata} (On Delisle's geographic projection used in the general map of the Russian empire)  \cite{Euler-pro-Desli-1777}; 
\item  \emph{De repraesentatione superficiei sphaericae super plano} (On the representation of  spherical surfaces on a plane) \cite{Euler-rep-1777};
\item
  \emph{De proiectione geographica superficiei sphaericae} (On the geographical projections of spherical surfaces) \cite{Euler-pro-1777};
\end{enumerate}
We have analyzed the first memoir in the preceding section.

  In the second memoir \cite{Euler-rep-1777}, questions on geographical maps are included in the setting of differential calculus and the calculus of variations. Let us recall in this connection that the problem of finding the shortest lines on a surface, which is typically a problem of calculus of variations, was  considered by Euler in his early youth, see \cite{E9}. In many ways, this problem is related to questions posed by geographical maps.

  In the introduction of his memoir \cite{Euler-rep-1777}, Euler declares that he shall not consider only maps obtained by central projection of the sphere onto the plane, which he calls ``optical projections", but mappings ``in the widest sense of the word."  He writes (\S 1 of \cite{Euler-rep-1777}):
 ``I take the word `mapping' in the widest possible sense; any point of the spherical surface is represented on the plane by any desired rule, so that every point of the sphere corresponds to a specified point in the plane, and conversely."

 In \S 9, Euler proves that there is no ``perfect" or ``exact" mapping from the sphere onto a plane. The meaning of the word ``perfect" was the cause of some confusion among the authors who quoted this memoir by Euler (it was interpreted as the fact that the map cannot preserve distances up to a scale, a fact that was well-known long before Euler), see the report in the paper \cite{Geography-Galkin}.   A precise statement of Euler's theorem is given in the paper \cite{CP}. In this context, a map is said to be \emph{perfect} if, locally, distances are preserved infinitesimally along the meridians and the parallels, and if, furthermore, the angles between meridians and parallels are preserved. Thus, we see once again the importance of these two systems of lines.
   Euler proves this result through the study of certain partial differential equations. A detailed and modern proof based on Euler's ideas is  given in \cite{CP}. We stress on the fact that, although meridians are geodesics, parallels are not, therefore one cannot immediately conclude that the map considered by Euler is a local isometry.

 With this result established, Euler searches for maps which best approximate the desired properties. He examines several particular projections of the sphere, looking systematically for the partial differential equations that they satisfy. 
He considers in particular  the following three kinds of maps:

\begin{enumerate}
\item \label{map1} maps in which the images of all the meridians are perpendicular to a given axis (the ``horizontal" axis in the plane), while all parallels are parallel to it; 

\item \label{map3} conformal maps;

\item \label{map2} maps where surface area is preserved (up to scale).

\end{enumerate}
He then gives examples of various maps satisfying each of the above three properties:  projections of the sphere onto a tangent plane, onto a cylinder  tangent to the equator, etc. and he studies distance and angle distortion under them. At the end of his memoir (\S 60), he notes that his work is rather theoretical, and has no immediate practical use.

  In contrast, in the memoir \cite{Euler-pro-1777}, Euler studies projections that are used in practical geography, that is, in the actual construction of maps (he emphasizes this in \S 20). Unlike the memoir \cite{Euler-rep-1777} which is based on the techniques of differential calculus, the present  one makes heavy use of spherical trigonometry. Euler gives formulae for the images of the equator, parallels and meridians under stereographic projection onto a plane that is tangent to the sphere at an arbitrary point.

  Euler's duties as a geographer also involved the examination of thousands of maps that were printed by the Academy. He was displeased with this task and he pretended that it was responsible for the deterioration of his vision. In a letter to Christian Goldbach, dated August 21st (September 1st, new Style), 1740, he writes (cf.  \cite{Euler-Goldbach} p. 163, English translation p. 672):
\begin{quote}\small
Geography is fatal to me. As you know, Sir, I have lost an eye working on it; and just now I nearly risked the same thing again. This morning I was sent a lot of maps to examine, and at once I felt the repeated attacks. For as this work constrains one to survey a large area at the same time, it affects the eyesight much more violently than simple reading or writing. I therefore most humbly request you, Sir, to be so good as to persuade the President by a forceful intervention that I should be graciously exempted from this chore, which not only keeps me from my ordinary tasks, but also may easily disable me once and for all.  \end{quote}

 In the next section, we shall talk about the work of Lagrange on geography.
 
 \section{Lagrange}\label{s:Lagrange}
    
   Joseph-Louis  Lagrange was 19 years younger than Euler and he was probably, among the large number of mathematicians with whom the latter corresponded, the only one capable of  understanding all his works. The two men, in their long mathematical correspondence, discussed a great variety of topics, including number theory, analysis, geometry, physics, etc. In almost every field, they obtained complementary results, and geography was no exception.  
   
 Two years after the three memoirs of Euler that we mentioned  appeared in print, Lagrange published two memoirs on geography, \cite{Lagrange1779}, both
 carrying the title
\emph{Sur la construction des cartes g\'eographiques} (On the construction of geographical maps).
In these memoirs, Lagrange declares that he extends works of several mathematicians, and he mentions in particular Euler and Lambert.\footnote{Johann-Heinrich Lambert (1728-1777)  was a Swiss-born mathematician for whom Euler had a great consideration. He was completely self-taught and he became a remarkable mathematician, astronomer, physicist and philosopher. He was hired at the Berlin Academy of Sciences at Euler's recommendation, and he spent there the last ten years of his life. Lambert 
  is sometimes considered as the founder of modern cartography.  His work \emph{Anmerkungen und Zus\"atze zur Entwerfung der Land- und Himmelscharten} (Remarks and complements for the design of terrestrial and celestial maps, 1772) \cite{Lamb-Anmer} contains seven new types of geographical maps, each one having important features. His name is associated with the so-called Lambert conformal conic projection, the transverse Mercator, the Lambert azimuthal equal area projection, and the Lambert cylindrical equal-area projection. In the same memoir, Lambert obtained a mathematical characterization of an arbitrary angle-preserving map from the sphere onto the plane. 
His projections are still mentioned in the modern textbooks of cartography, and some of them were still in use in the twentieth century, for military and other purposes, until the appearance of satellite maps. Lambert, in his work, took into account the fact that the Earth is spheroidal and not spherical. His memoir  \cite{Lamb-Anmer} is part of his larger treatise \emph{Beitr\"age zum Gebrauche der Mathematik und deren Anwendung durch} (Contributions to the use of mathematics and its applications) \cite{Lambert-Bey}.}

Besides extending the works of his predecessors in geography, Lagrange brought several  new ideas of his own. In particular, he introduced a notion of  distortion that had an important impact on the works of several mathematicians, in particular on Pafnouti Chebyshev, who worked in Saint Petersburg a few decades after Euler.  We shall mention Chebyshev's contribution to geography in the concluding section of this paper. Let us now review some works of Lagrange related to our topic.

 In the introduction to his first memoir  \cite[p. 637]{Lagrange1779}, Lagrange starts by declaring that a geographical map is nothing but a plane figure which represents the surface of the Earth, or some part of it. He recalls that since the Earth is spherical---or rather spheroidal, as he says---, it is not possible to represent on a Euclidean plane an arbitrary part of it without altering the positions and distances of various places. One can look then for best possible maps; and the less the alteration of the distances is, the better the geographical map is. He then surveys several projections of the Earth (and of the Celestial sphere) that were used before him, among which we have the \emph{central projection},  whose center of projection is at the center of the globe. This projection sends the great circles (therefore, the meridians) to straight lines, whereas the small circles (such as the parallels) are sent either to circles or to ellipses, according to whether their plane is parallel or not to the projection plane. A major advantage of the central projection is that the shortest way between two locations is the straight line on the map, since the great circles of the sphere are sent to straight lines. Lagrange recalls that this projection is usually applied in such a way that the projection plane is parallel to the equator, and in this case parallels are sent to circles. This projection, he says, is mostly used for maps of the Celestial sphere.

Lagrange then mentions the stereographic projection from a point on the sphere onto a plane, and he recalls two of its main properties: it sends circles to circles, and it is conformal (that is, it preserves angles).  He attributes this projection to Ptolemy \cite[p. 639]{Lagrange1779}.\footnote{According to Delambre \cite{Del}, Vol. 1, p.
184ff., d'Avezac \cite{Dav}, p. 465 and Neugebauer \cite{Neu, Neu1948}, p. 246., the stereographic projection from the sphere onto a plane where the center is taken to be
a point on the sphere (say the North pole) and where the plane passes through the
center and is perpendicular to the radius passing through the center of projection,  which was probably the
most popular projection of the sphere among mathematicians, was already used
by Hipparchus back in the second century BCE.} He declares that the latter was aware of the first of these properties, which he describes in his treatise  \emph{Sphaerae a planetis projectio in planum},\footnote{This is the work  known as the \emph{Planisphaerum} which we mentioned in \S \ref{s:antiquity}.}  but that the angle-preserving property may not have been noticed by the Greek astronomer. 
 
 After introducing mappings from the sphere onto a plane that are projections from some center (allowing also the center to be at infinity, in which case he calls the projection \emph{orthographic}), Lagrange  mentions other conformal projections that are not stereographic. He notes that there are infinitely many such projections, and he also considers much more general maps. He declares, like Euler did before him (cf. the introduction of his memoir \cite{Euler-rep-1777}), that geographical maps may be arbitrary maps from (some part of) the sphere onto a plane. This leads him to the  following question:  
 \medskip
 
 \centerline{\emph{What is a
  general mapping between two surfaces?}}
  \medskip
  
    We mention incidentally that this important question is a two-dimensional analogue of the question of \emph{what is a general function} defined on the real line (or an interval) and which, at that time, was the subject of a fierce debate among mathematicians; see the account in \cite{Riemann0}.

     Returning to maps from the sphere onto the Euclidean plane, Lagrange writes   \cite[p. 640]{Lagrange1779}  that the only thing we require for drawing a map is to specify the images of meridians and parallels according to a certain rule. Then, we may plot the various places relatively to these lines, in the same way this is done on the surface of the Earth with respect to the circles of longitude and latitude. Now, the images of the meridians and parallels are no more restricted to be circles or lines. They can be, using Lagrange's terms,  ``mechanical lines," that is, lines drawn by any mechanical device. In the classical terminology used by the Greeks for curves, this means that they can be arbitrary lines. He then discusses the reduced marine maps in which the images of the meridians and of the parallels are parallel straight lines, and where the ratios of the degrees of latitude and longitude on these images are  the same as on the sphere.

 Lagrange then makes a short historical overview of the question of how general can be a mapping from the sphere onto the plane that may be used in map drawing. He refers to the work of Lambert, who was the first to talk about arbitrary angle-preserving maps from the sphere to the plane in relation with geography. He recalls that the latter expressed the idea of the determination of the images of the meridians and the parallels by the sole condition that the map is angle-preserving. In fact, in
   his memoir \cite{Lambert-Bey}, Lambert solved the problem of characterizing least-distortion maps among those which are angle-preserving. Lagrange also recalls that Euler, like Lambert did before him,  gave a solution of the problem of finding maps with least distortion among arbitrary angle-preserving maps. Lagrange then develops his own solution, by a method which is different from those of Lambert and Euler. He asserts that none of his predecessors has considered yet the general problem of determining all the conformal maps by which the images of the meridians and the parallels are circles and he considers this problem in some detail. As a matter of fact, Lagrange, in his paper, solves the problem of finding all the orthonormal projections of a surface of revolution which send meridians and parallels to straight lines or circles. He introduces an explicit formula for the local distortion factor of  a conformal map, as the ratio of the infinitesimal length element at the image by the infinitesimal length element at the source. This formula played a central role in the work of Chebyshev, as we shall recall in the next section.

\section{A glimpse of later developments}\label{s:Conclusion}

In the preceding sections, we saw that in the eighteenth century, cartography underwent substantial developments, thanks to the efforts of prominent mathematicians like Euler, Lambert and Lagrange. The subject  continued to grow in the nineteenth century. In this concluding section, I would like to say a few words about these developments.

We start with the work of Gauss.

Gauss considered himself as a physicist rather than a mathematician. We recall in this respect that the terminology ``isothermal coordinates" which he introduced in the study of the differential geometry of surfaces and which is still used today, 
for a locally conformal map between a
subdomain of the plane and a subdomain of the surface, clearly indicates the fact that while working on this subject, he was thinking about heat diffusion.
Gauss was also 
in charge of surveying geodetically the German kingdom of
Hannover.  
In 1825, he published a paper, in the \emph{Astronomische Abhandlungen}  (Memoirs on astronomy), titled \emph{Allgemeine L\"osung der Aufgabe, die Teile einer gegebenen Fl\"ache auf eine ander gegebene Fl\"ache so abzubilden dass die Abbildung dem Abgebildeten in den kleinisten Teilen \"ahnlich wird.} (General solution of the problem: to represent the parts of a given surface on another so that  the smallest parts of the representation shall be similar to the corresponding parts of the surface represented) \cite{Gauss-Copenhagen}. The title indicates enough explicitly the relation with geography. As we saw, this is the kind of problem that was considered by Euler, Lagrange and other mathematicians, in their work on this field. 
   In this paper, Gauss showed that every sufficiently small neighborhood  of a point in an arbitrary real-analytic surface can be mapped conformally onto a subset of the plane,  however, without solving the problem of mapping conformally an arbitrary finite portion of the surface;
this was one of the questions considered by Riemann, after Gauss. In the preface, Gauss wrote that his aim was only to construct geographical maps and to study the
general principles of geodesy for the task of land surveying. 

Surveying the kingdom
of Hannover took nearly two decades to be completed, and it led Gauss gradually to the
investigation of triangulations, to the use of the method of least squares in geodesy, 
and then to his famous memoir \emph{Disquisitiones generales circa superficies curvas} \cite{Gauss-English}, in which we can find coordinates that he computed of several cities in Germany (\S 27; p. 43 of the English translation \cite{Gauss-English}). 

Let us mention two other papers by Gauss related to geodesy:  \emph{Bestimmung des Breitenunterschiedes zwischen den Sternwarten von G\"ottingen und Altona durch Beobachtungen am Ramsdenschen Zenithsector} (Determination of the latitudinal difference between
the observatories in G\"ottingen and Altona by observations with a
Ramsden zenith sector) \cite{Gauss-Bestimmung} (1928) and \emph{Untersuchungen \"uber Gegenst\"ande der h\"ohern Geod\"asie} (Research on objects of higher geodesy) \cite{Gauss-Untersuchungen} (1843 and 1847). In the latter, Gauss uses the method of least squares.

Another major representative of the 19th-century differential-geometric work on geography is Eugenio Beltrami, who worked on this subject in the tradition of Gauss. His 1865 paper \emph{Risoluzione del problema: Riportare i punti di una superficie sopra un piano in modo che le linee geodetiche vengano rappresentate da linee rette} (Solution of the problem: to send the points of a surface onto a plane in such a way that the geodesic lines are represented by straight lines) \cite{Beltrami1865} contains his well-known result saying that a Riemannian metric on a surface that can be locally mapped onto the plane in such a way that the geodesics are sent to Euclidean lines has necessarily constant curvature. Beltrami declares  that a certain amount of the research done before him on this kind of problem was directed towards questions of conservation of angles or of area, and that even though these two properties are considered to be the simplest and most important ones for geographical maps, there are other properties that one might want to preserve. 
 He mentions that the central projection of the sphere is the only map that transforms the geodesics of the sphere into Euclidean straight lines, but that one needs other maps, in  which the images of great circles are not very far from being straight lines.  He writes that beyond its applications to geographical maps, this kind of problem will lead to  ``a new method of geodesic calculus, in which the questions concerning geodesic triangles on surfaces can all be reduced to simple questions of plane trigonometry."

We now say a few words on the work of Chebyshev, one of those prominent 19th century mathematicians who worked in Saint Petersburg and who is considered as the  founder of the Russian school there, which in some sense replaced the school founded by Euler. 

  Chebyshev was interested in applications of mathematics, and geography was among them. 
  In 1856, he wrote two papers on geography \cite{Cheb1, Cheb2} carrying the same title as the two papers that Lagrange  published 77 years before \cite{Lagrange1779}, \emph{Sur la construction des cartes g\'eographiques}.  In these papers, Chebyshev, addressed the problem of finding geographical maps whose distortion is minimal. He made a relation between this problem and  Laplace's equation, thus reducing the problem of finding the best geographical map to a problem in potential theory. 
At the beginning of the second paper, Chebyshev declares: ``Today, [Mathematical sciences] produce a greater interest because of their influence on art and industry. Not only practice makes profit of these relations: conversely, science itself grows under the influence of practice."  Then, Chebyshev elaborates on the importance in this context of the problem of constructing geographical maps, establishing relations with the problem of heat distribution and with other problems and including them in the same setting of infinitesimal calculus.

In his papers on geography, Chebyshev, starting with formulae that Lagrange gave on the distortion of maps from the sphere to the plane, announced several important results, including the existence of a map with least distorsion from an arbitrary simply connected open subset of the sphere bounded by a twice differentiable curve, onto to the Euclidean plane. This map has the property that the parallels and the meridians are represented by curves that do not differ significantly from circles or straight lines.
Chebyshev also proved that such a map is  unique up to a similarity of the Euclidean plane. Furthermore, he showed that the magnification ratio of this  map is constant on the boundary curve of the domain and he obtained precise estimates on the distorsion in the case where the boundary curve is close to an ellipse. Chebyshev did not publish the proofs of his results.
     
     Darboux  in a paper published in 1911 \cite{Darboux-Chebyshev} and which carries the same title as Chebyshev's, \emph{Sur la construction des cartes g\'eographiques}, gave complete proofs of Chebyshev's main result. The same result was also reviewed in a paper \cite{Milnor} by Milnor, where Chebyshev's work, together with that of Lagrange on geography, are put in modern  perspective. We refer the interested reader to Chebyshev's papers \cite{Cheb1, Cheb2} and the reviews of these papers in \cite{2016-Tchebyshev} and \cite{Papa-qc}.

Talking about the 19th century, let us also mention the second doctoral thesis of the French mathematician Ossian Bonnet,  titled \emph{Sur la th\'eorie math\'ematique des cartes g\'eographiques} (On on the mathematical theory of geographical maps) \cite{Bonnet-these}. In this work, Bonnet starts by recalling that the first geographical maps were projections that are subject to the laws of perspective, and he then declares that ``some astronomers" abandoned this perspective way of drawing maps, and considered arbitrary ways of drawing the lines corresponding to the meridians and parallels, depending on the usage of the map. He attributes to Lambert the formulation of the problem of finding, for a given region of the sphere, the images of the meridians and the parallels by a map which is conformal and such that at the infinitesimal level distances are preserved on these lines. Bonnet declares that Lambert did not completely solve  the problem, that Euler and Lagrange considered it again (we talked about a similar problem considered by Euler in the paper  \cite{Euler-rep-1777}), and that Gauss, in his \emph{M\'emoire Couronn\'e} \cite{Gauss-s}, solved it in full generality, without making any hypothesis on the form of the Earth. The main part of Bonnet's thesis contains an exposition of the work of Lambert on cartography, together with a simplification of the solutions of the questions he asked that were given by Lagrange, Euler and Gauss.
     
       Bonnet's  thesis was published in the \emph{Journal de math\'ematiques pures et appliqu\'ees}, edited by Liouville.
In a note at the end of the paper,  the latter writes  that he addressed these questions on geography in a series of lectures he
gave at the Coll\`ege de France, in the academic year 1850-1851, and that he
wishes to publish them. He also notes that he presented his ideas on
the subject in the Notes of his edition of Monge's book, \emph{Applications de l'analyse \`a la g\'eom\'etrie}, \cite{Monge}. In fact, in Note V and VI, Liouville formulates a problem he  calls the 
 ``three-dimensional geographical
drawing problem". The title of Note V is \emph{Du trac\'e g\'eographique des surfaces les unes sur
les autres} (On the geographical drawing of surfaces one onto the other). In this note, Liouville formulates the problem as the one of
finding a mapping between two surfaces which is a similarity at the infinitesimal
level. 
  He declares that this is equivalent to requiring that the infinitesimal triangles on
the first surface are sent by the map to similar infinitesimal triangles on the second
one. He then formulates the same problem using ratios of  infinitesimal line elements at a point on one surface and the one at its image: the ratio at any point should not depend
on the chosen direction. He notes that this condition (which, in fact, is conformality) was adopted by
 Lambert, Lagrange and Gauss as a general principle in their theory of geographical maps. In Notes V and VI, Liouville gives his solution to the problem he formulates.

  As a conclusion, we quote Darboux, from his his talk at the 1908 ICM held in Rome, titled \emph{Les origines, les m\'ethodes et les probl\`emes de la g\'eom\'etrie infinit\'esimale} (The origins, methods and problems of infinitesimal geometry) \cite{Darboux-ICM}: \begin{quote} \small 
Like many other branches of human knowledge, infinitesimal geometry was born in the study of practical problems. The Ancients were already busy in obtaining plane representations of the various parts of the Earth, and they had adopted the idea, which was so natural, of projecting onto a plane the surface of our globe. During a very long period of time, people were exclusively attached to these methods of projection, restricting simply to the study of the best ways to choose, in each case, the point of view and the plane of the projection. It was one of the most penetrating geometers, Lambert, the very estimated colleague of Lagrange at the Berlin Academy, who, pointing out for the first time a property which is common to the Mercator maps, also called \emph{reduced maps}, and to those which are provided by the stereographic projection, was the first to conceive the theory of geographical maps from a really general point of view. He proposed, with all its scope, the problem of representing the surface of the Earth on a plane keeping the similarity of the infinitely small elements. This beautiful question, which gave rise to the researches of Lambert himself, of Euler, and to two very important memoirs of Lagrange, was treated for the first time in all generality by Gauss. [...] Among the essential notions introduced by Gauss, one has to note the systematic use of the curvilinear coordinates on a surface, the idea of considering a surface like a flexible and inextensible fabric, which led the great geometer to his celebrated theorem on the invariance of total curvature, to the beautiful properties of geodesic lines and their orthogonal trajectories, to the generalization of the theorem of Albert Girard on the area of the spherical triangle, to all these concrete and final truths which, like many other results due to the genius of the great geometer, were meant to preserve, across the ages, the name and the memory of the one who was the first to discover them.
\end{quote} 

In many ways, this article may be considered as an expansion of this quote.
\bigskip

\noindent {\bf Acknowledgements.} This paper was written after a talk I gave on the same subject at the International Conference on History and Recent Developments in Mathematics  that took place on December 17--19 at the Madhuben and Bhanubhai Patel Institute of Technology (New VV Nagar, Gujarat, India). I would like to thank Darshana Prajapati and S. G. Dani who organized this event. I would also like to thank the referee of this paper for his very useful comments.

\end{document}